\documentclass{article}

\usepackage{BICA}
\usepackage{pgf,tikz}
\usepackage{verbatim}
\usetikzlibrary{decorations.pathreplacing,calligraphy}
\usepackage{longtable}

 \newtheorem{Theorem}{Theorem}[section]
   
   \newtheorem{Corollary}[Theorem]{Corollary}

   \newcommand{\qed}{\hphantom{.}\hfill \rule[1pt]{8pt}{8pt}\medbreak}

\begin{document}

\title{ A 2-Person Game Decomposing 2-Manifolds}
\author{David R. Berman}
\affil{Wilmington, NC
\newline
\email{ bermand@uncw.edu} }
\author{Lee O. Leonard Jr}
\affil{Austin, TX
\newline
\email{lleonardjr@icloud.com}}
\date{June 30, 2024}

\maketitle

\begin{abstract}
Two players play a game by alternately splitting a surface of a compact $2$-manifold along a simple closed curve 
that is not null-homotopic and attaching disks to the resulting boundary; the last player who can move wins. 
Starting from an orientable surface, the $G$-series is $01\dot{2}\dot{0}$ according to increasing genus. 
Starting from a nonorientable surface, the $G$-series is $012\dot{4}60\dot{3}$ according to increasing genus. 
Nim addition determines the $G$-values of the remaining compact $2$-manifolds.
\end{abstract}


\section{Introduction}

In this paper we introduce a two player game played on a two-dimensional topological manifold, and give a winning strategy. 
We begin by reviewing necessary material from combinatorial game theory and topology. 

Combinatorial games have a rich history in recreational mathematics and number theory. The definitive book is Winning Ways for Your Mathematical Plays by Berlekamp, Conway, and Guy~\cite{BCG},
which contains theory, history, and numerous games that are played and analyzed. Our game involves two players who move alternately, the last one who can move wins. The game terminates
after a finite number of moves. The game is impartial, meaning the same moves are available to each player. The quintessential such game is \textit{Nim}.

Nim is played with a finite set of heaps, each containing a finite number of counters. On a move, a player chooses a heap and removes one or more counters. The last player who can move wins.
A winning strategy is found by calculating the nim sum of the heaps. The \textit{nim sum} of  $a$ and $b$ is denoted $a \oplus b$ and is found by writing the numbers in binary and adding 
with no carries. Second player wins if the nim sum is zero while first player wins if the sum is positive. Suppose the sum is positive. Then find the largest power 
of two that is unpaired. Write as 1 + all smaller powers of two. First player removes powers that are already paired, plus 1. Now the powers are all paired and the nim sum is zero. So the first 
player, now as second player, will win. If the sum is zero, all powers are paired, so removing counters from a single heap will leave some power unpaired, and the nim sum will be positive. For 
example with heap sizes 21, 11, and 6, $21 \oplus 11 \oplus 6 = (1 + 4 + 16) \oplus (1 +2 + 8) \oplus (2 +4) = 8 + 16 = 24$, and first player wins by removing $8$ from the heap $21$.  

A winning strategy for any two player combinatorial impartial terminating game can be determined using the Sprague-Grundy theory, which assigns a non-negative integer, called the \textit{Grundy number}, to each game position. 
For short we will call the Grundy number the value. The value of terminal positions, that is with no moves, is 0. These are losses for the current player. Proceed inductively. 
The value of a given game position is the minimal excluded integer, i.e. $mex$, of the value of the game positions after each possible move. For example, $mex\{1,3,5\} = 0$, and $mex\{0,1,3,5\} = 2$.
Finally, if a game is the union of subgames then the value of the game is the nim sum of the values of the subgames, and the analysis is similar to nim. 
A winning strategy can now be stated. If the initial game value is positive, then the first player wins by moving to a position with value 0. If the initial game value is 0, then the second player wins
because every move results in a position of positive value. 


Classification of compact surfaces, i.e., compact and connected $2$-manifolds, is well known and is stated below. We follow the exposition given in Massey~\cite{M1}, and have used~\cite{G, HW, M2, V} as well. 
A compact $2$-manifold is composed of a finite disjoint collection of compact surfaces. The \textit{connected sum} of two surfaces is formed by removing a disk from each surface and joining the surfaces along the boundary.
The basic compact surfaces are the \textit{sphere}, \textit{torus}, and \textit{projective plane}. The sphere and torus are orientable and the projective plane nonorientable. We also use the \textit{annulus} and \textit{M\"obius band} 
which are bordered surfaces. Standard constructions of these are shown below. Edges of each square with the same label are identified in the direction of the arrow. 
\\[.05 in]

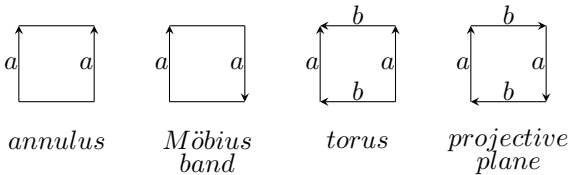
\begin{tikzpicture}


\draw [-stealth] (1,1)--(1,2);\draw (1,2) -- (2,2);\draw [stealth-] (2,2)--(2,1);\draw (2,1)--(1,1);
\node at (1.5,.5){$annulus$};\node at (.9,1.5){$a$};\node at (1.9,1.5){$a$};

\draw [-stealth] (3,1)--(3,2);\draw (3,2) -- (4,2);\draw [-stealth] (4,2)--(4,1);\draw (4,1)--(3,1);
\node at (3.5,.5){$M\ddot{o}bius$};\node at (3.5,.2) {$band$};\node at (2.9,1.5){$a$};\node at (3.9,1.5){$a$};

\draw [-stealth] (5,1)--(5,2);\draw [stealth-](5,2) -- (6,2);\draw [stealth-] (6,2)--(6,1);\draw [-stealth](6,1)--(5,1);
\node at (5.5,.5){$torus$};\node at (4.9,1.5){$a$};\node at (5.9,1.5){$a$};\node at (5.5,1.15){$b$};\node at (5.5,2.15){$b$};

\draw [-stealth] (7,1)--(7,2);\draw [-stealth](7,2) -- (8,2);\draw [-stealth] (8,2)--(8,1);\draw [-stealth](8,1)--(7,1);\node at (7.5,.5){$projective$};
\node at (7.5,.2){$plane$};\node at (6.9,1.5){$a$};\node at (7.9,1.5){$a$};\node at (7.5,1.15){$b$};\node at (7.5,2.15){$b$};

\end{tikzpicture}
\\[.05 in]
We will need the observation that a \textit{M\"obius band} with a disk attached to the edge is homeomorphic to the projective plane.
A beautiful pictorial proof is given in~\cite{DL}. The disk is called a \textit{cross-cap}. 

\begin{Theorem}\label{classification}
Every compact surface is homeomorphic to either a sphere, a connected sum of $g$ tori, or a connected sum of $g$ projective planes, where $g$ is the genus. The genus of a sphere is 0.  
\end{Theorem}

We denote the connected sum of $g$ tori by $\mathit{og}$, the orientable surface of genus $g$, and the connected sum of $g$ projective planes by $\mathit{ng}$, the nonorientable surface of genus $g$. 
For convenience we write $o0 = n0$ for the sphere. The genus is related to the Euler characteristic which is not needed in our presentation. The proof of the classification theorem in Massey uses the following simplification. 

\begin{Theorem}
The connected sum of a torus and a projective plane is homeomorphic to the connected sum of three projective planes. 
\end{Theorem}

We will make use of the following corollary.
\begin{Corollary}
If  $g=a+b, b>0$, and $a$ is even, then the connected sum of $g$ projective planes is homeomorphic to the connected sum of $a/2$ tori and $b$ projective planes. 
\end{Corollary}

We are now ready to introduce our whimsical game. While not playable in practice due to the difficulty in visualizing nonorientable surfaces, it is nevertheless interesting to study and analyze. 
\\[.05 in]

\section{Manifold Decomposition Game}

Assume a 2-manifold is given and is composed of compact surfaces. 
Two players play a game by alternately forming a proper decomposition of a compact surface. A proper decomposition starts with an essential simple closed curve, denoted $J$. 
Essential means the curve is not null-homotopic, or equivalently does not bound a disk in the manifold (see Leonard~\cite{L}). 
The simple closed curve is expanded into a tubular, i.e., regular neighborhood. The neighborhood is either an annulus or a M\"obius band. Remove the interior of the 
neighborhood leaving the boundary, two circles for the annulus and one circle for the M\"obius band. The decomposition, i.e., move, is completed by capping off the circles. 
Each move results in either a single surface or a pair of surfaces.
The game begins with either $og$, the orientable surface of genus $g$, or $ng$, the nonorientable surface of genus $g$. The game ends with a collection of spheres, 
because there are no essential simple closed curves on a sphere. 

\begin{Theorem}
For genus $g$, the following decompositions are always available and comprise all possible proper decompositions of a compact surface. These are the moves in the manifold decomposition game. 

\begin{tabular}{lll}
a) &  $og \rightarrow o(g-1)$ 	      &      $g - 1 \geq 0$  \\
b) &  $og \rightarrow (oa, ob)$               &      $a, b >0,  a+b=g$ \\
c) &  $ng \rightarrow n(g-1)$                &      $g-1 \geq 0$  \\
d) &  $ng \rightarrow o((g-1)/2)$          &       $g-1 \geq 0$ and even  \\
e)  & $ng \rightarrow n(g-2)$                &       $g-2 \geq 0$  \\
f)  & $ng \rightarrow o((g-2)/2)$            &       $g-2 \geq 0$ and even  \\
g)  &  $ng \rightarrow (na, nb)$               &      $a, b>0,  a+b=g$ \\
h)  & $ng \rightarrow (o(a/2), nb)$          &       $a, b > 0, a$ even,  $a+b=g$ \\
\end{tabular}

\end{Theorem}
\begin{Proof} 
The cases are visualized below. Cases $d, f,$ and $h$ make use of $e$ and Corollary 1.3.
In cases $c$ and $d$, the curve $J$ is identified antipodally. The tubular neighborhood is a M\"obius band. The projective plane, i.e., crosscap, is composed
of the M\"obius band and the boundary. When the M\"obius band is removed and the boundary capped, the projective plane disappears.
In cases $e$ and $f$, the curve $J$ involves two crosscaps. The tubular neighborhood joins the two M\"obius bands forming an annulus and capping the boundaries with disks
removes both projective planes.  
\\[.05 in]
$a)$

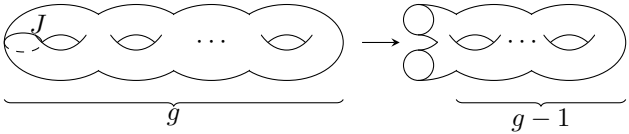
\begin{tikzpicture}[scale=1/2]

\begin{scope}\clip (-2,-2) rectangle (1,2);\draw (0,0) ellipse (1.5 and 1);\end{scope}
 \node at (-.6,.5) {$J$};
\draw (-1.5,0) arc (180:0:.5 and .25); \draw[dashed] (-1.5,0) arc (180:360:.5 and .25);
\begin{scope}\clip (0,-.9) ellipse (1 and 1.5);\draw (0,1.3) ellipse (1 and 1.5);\end{scope}
\begin{scope}\clip (0,1.3) ellipse (1 and 1.5);\draw (0,-1.3) ellipse (1 and 1.5);\end{scope}

\begin{scope}\clip (1,-2) rectangle (3,2);\draw (2,0) ellipse (1.5 and 1);\end{scope}
\begin{scope}\clip (2,-.9) ellipse (1 and 1.5);\draw (2,1.3) ellipse (1 and 1.5);\end{scope}
\begin{scope}\clip (2,1.3) ellipse (1 and 1.5);\draw (2,-1.3) ellipse (1 and 1.5);\end{scope}\draw (4,0) node {$\cdots$};

\begin{scope}\clip (3,-2) rectangle (5,2);\draw (4,0) ellipse (1.5 and 1);\end{scope}

\begin{scope}\clip (5,-2) rectangle ((8,2);\draw (6,0) ellipse (1.5 and 1);\end{scope}
\begin{scope}\clip (6,-.9) ellipse (1 and 1.5);\draw (6,1.3) ellipse (1 and 1.5);\end{scope}
\begin{scope}\clip (6,1.3) ellipse (1 and 1.5);\draw (6,-1.3) ellipse (1 and 1.5);\end{scope}

\begin{scope} \clip (9.5,-2) rectangle (11,2);
\begin{scope}\clip (8,-2) rectangle (10.5,2);\draw (9.5,0) ellipse (1.5 and 1);\end{scope}
\begin{scope} \clip (9.5, -1) rectangle (10,.3);\begin{scope}\clip (9.5,-.9) ellipse (1 and 1.5);\draw (9.5,1.3) ellipse (1 and 1.5);\end{scope}
\begin{scope}\clip (9.5,1.3) ellipse (1 and 1.5);\draw (9.5,-1.3) ellipse (1 and 1.5);\end{scope}\end{scope}
\end{scope}
\begin{scope}\clip (11,-.9) ellipse (1 and 1.5);\draw (11,1.3) ellipse (1 and 1.5);\end{scope}
\begin{scope}\clip (11,1.3) ellipse (1 and 1.5);\draw (11,-1.3) ellipse (1 and 1.5);\end{scope}
\draw (9.5,.6) circle (.4); \draw (9.5,-.6) circle (.4);

\begin{scope}\clip (10.5,-2) rectangle (12.5,2);\draw (11.5,0) ellipse (1.5 and 1);\end{scope}

\begin{scope}\clip (12.5,-2) rectangle ((15.5,2);\draw (13.5,0) ellipse (1.5 and 1);\end{scope}
\begin{scope}\clip (13.5,-.9) ellipse (1 and 1.5);\draw (13.5,1.3) ellipse (1 and 1.5);\end{scope}
\begin{scope}\clip (13.5,1.3) ellipse (1 and 1.5);\draw (13.5,-1.3) ellipse (1 and 1.5);\end{scope}\draw (12.25,0) node {$\cdots$};

\draw[-stealth] (8,0) -- (9,0);
\draw [decorate, decoration = {brace,mirror}] (-1.5,-1.5) -- (7.5,-1.5);\draw (3,-1.5) node[below] {$g$};
\draw [decorate, decoration = {brace,mirror}] (10.5,-1.5) -- (15,-1.5);\draw (12.75,-1.5) node[below] {$g-1$};

\end{tikzpicture}
\\[.05 in]
\indent $b)$

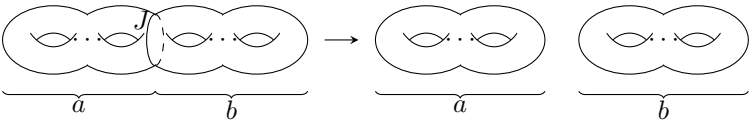
\begin{tikzpicture}[scale=.45]
\begin{scope}\clip (-2,-2) rectangle (1,2);\draw (0,0) ellipse (1.5 and 1);\end{scope}
\begin{scope}\clip (0,-.9) ellipse (1 and 1.5);\draw (0,1.3) ellipse (1 and 1.5);\end{scope}
\begin{scope}\clip (0,1.3) ellipse (1 and 1.5);\draw (0,-1.3) ellipse (1 and 1.5);\end{scope}
 
\begin{scope}\clip (1,-2) rectangle (3,2);\draw (2,0) ellipse (1.5 and 1);\end{scope}
\begin{scope}\clip (2,-.9) ellipse (1 and 1.5);\draw (2,1.3) ellipse (1 and 1.5);\end{scope}
\begin{scope}\clip (2,1.3) ellipse (1 and 1.5);\draw (2,-1.3) ellipse (1 and 1.5);\end{scope}\draw (1,0) node {$\cdots$};

\begin{scope}\clip (3,-2) rectangle (5,2);\draw (4,0) ellipse (1.5 and 1);\end{scope}
\begin{scope}\clip (4,-.9) ellipse (1 and 1.5);\draw (4,1.3) ellipse (1 and 1.5);\end{scope}
\begin{scope}\clip (4,1.3) ellipse (1 and 1.5);\draw (4,-1.3) ellipse (1 and 1.5);\end{scope}

\begin{scope}\clip (5,-2) rectangle ((8,2);\draw (6,0) ellipse (1.5 and 1);\end{scope}
\begin{scope}\clip (6,-.9) ellipse (1 and 1.5);\draw (6,1.3) ellipse (1 and 1.5);\end{scope}
\begin{scope}\clip (6,1.3) ellipse (1 and 1.5);\draw (6,-1.3) ellipse (1 and 1.5);\end{scope}\draw (5,0) node {$\cdots$};

\draw [dashed] (3,.75) arc (90:-90:0.25 and .75);
\draw (3,.75) arc (90:270: 0.25 and .75); \node at (2.6,.6) {$J$};

\begin{scope}\clip (9,-2) rectangle (12,2);\draw (11,0) ellipse (1.5 and 1);\end{scope}
\begin{scope}\clip (11,-.9) ellipse (1 and 1.5);\draw (11,1.3) ellipse (1 and 1.5);\end{scope}
\begin{scope}\clip (11,1.3) ellipse (1 and 1.5);\draw (11,-1.3) ellipse (1 and 1.5);\end{scope}
 
\begin{scope}\clip (12,-2) rectangle (15,2);\draw (13,0) ellipse (1.5 and 1);\end{scope}
\begin{scope}\clip (13,-.9) ellipse (1 and 1.5);\draw (13,1.3) ellipse (1 and 1.5);\end{scope}
\begin{scope}\clip (13,1.3) ellipse (1 and 1.5);\draw (13,-1.3) ellipse (1 and 1.5);\end{scope}\draw (12,0) node {$\cdots$};

\begin{scope}\clip (15,-2) rectangle (18,2);\draw (17,0) ellipse (1.5 and 1);\end{scope}
\begin{scope}\clip (17,-.9) ellipse (1 and 1.5);\draw (17,1.3) ellipse (1 and 1.5);\end{scope}
\begin{scope}\clip (17,1.3) ellipse (1 and 1.5);\draw (17,-1.3) ellipse (1 and 1.5);\end{scope}
 
\begin{scope}\clip (18,-2) rectangle (21,2);\draw (19,0) ellipse (1.5 and 1);\end{scope}
\begin{scope}\clip (19,-.9) ellipse (1 and 1.5);\draw (19,1.3) ellipse (1 and 1.5);\end{scope}
\begin{scope}\clip (19,1.3) ellipse (1 and 1.5);\draw (19,-1.3) ellipse (1 and 1.5);\end{scope}\draw (18,0) node {$\cdots$};

\draw[-stealth] (8,0) -- (9,0);
\draw [decorate, decoration = {brace,mirror}] (-1.5,-1.5) -- (3,-1.5);\draw (.75,-1.5) node[below] {$a$};
\draw [decorate, decoration = {brace,mirror}] (3,-1.5) -- (7.5,-1.5);\draw (5.25,-1.5) node[below] {$b$};
\draw [decorate, decoration = {brace,mirror}] (9.5,-1.5) -- (14.5,-1.5);\draw (12,-1.5) node[below] {$a$};
\draw [decorate, decoration = {brace,mirror}] (15.5,-1.5) -- (20.5,-1.5);\draw (18,-1.5) node[below] {$b$};

\end{tikzpicture}
\\[.05 in]
\indent $c)$

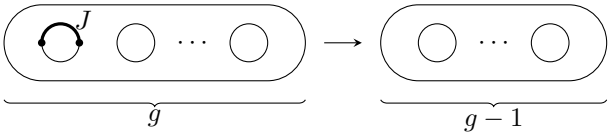
\begin{tikzpicture}[scale=1/2]

\draw (-1.5,1) -- (4.5,1) arc (90:-90:1) -- (-1.5,-1) arc (270:90:1);
\draw (-1,0) circle (.5); \draw (1,0) circle (.5); \draw (2.5,0) node {$\cdots$}; \draw (4,0) circle (.5);
\fill (-1.5,0) circle (.1); \fill (-.5,0) circle (.1); \draw[very thick] (-.5,0) arc (0:180:.5); \node at (-.4,.6) {$J$};

\draw (8.5,1) -- (12.5,1) arc (90:-90:1) -- (8.5,-1) arc (270:90:1);
\draw (9,0) circle (.5); \draw (12,0) circle (.5); \draw (10.5,0) node {$\cdots$}; 

\draw[-stealth] (6,0) -- (7,0);
\draw [decorate, decoration = {brace,mirror}] (-2.5,-1.5) -- (5.5,-1.5);\draw (1.5,-1.5) node[below] {$g$};
\draw [decorate, decoration = {brace,mirror}] (7.5,-1.5) -- (13.5,-1.5);\draw (10.5,-1.5) node[below] {$g-1$};
\end{tikzpicture}
\\[.05 in]
\indent $d)$

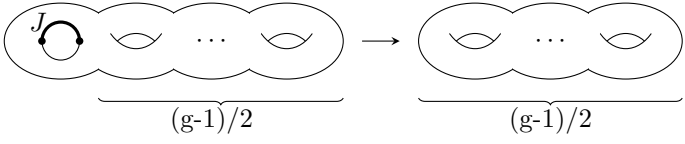
\begin{tikzpicture}[scale=1/2]
\begin{scope}\clip (-2,-2) rectangle (1,2);\draw (0,0) ellipse (1.5 and 1);\end{scope}
\draw (0,0) circle (.5);\draw[very thick] (.5,0) arc (0:180:.5); \node at (-.6,.5) {$J$};
\fill (-.5,0) circle (.1); \fill (.5,0) circle (.1); 

\begin{scope}\clip (1,-2) rectangle (3,2);\draw (2,0) ellipse (1.5 and 1);\end{scope}
\begin{scope}\clip (2,-.9) ellipse (1 and 1.5);\draw (2,1.3) ellipse (1 and 1.5);\end{scope}
\begin{scope}\clip (2,1.3) ellipse (1 and 1.5);\draw (2,-1.3) ellipse (1 and 1.5);\end{scope}\draw (4,0) node {$\cdots$};

\begin{scope}\clip (3,-2) rectangle (5,2);\draw (4,0) ellipse (1.5 and 1);\end{scope}

\begin{scope}\clip (5,-2) rectangle ((8,2);\draw (6,0) ellipse (1.5 and 1);\end{scope}
\begin{scope}\clip (6,-.9) ellipse (1 and 1.5);\draw (6,1.3) ellipse (1 and 1.5);\end{scope}
\begin{scope}\clip (6,1.3) ellipse (1 and 1.5);\draw (6,-1.3) ellipse (1 and 1.5);\end{scope}

\begin{scope}\clip (9.5,-2) rectangle (12,2);\draw (11,0) ellipse (1.5 and 1);\end{scope}
\begin{scope}\clip (11,-.9) ellipse (1 and 1.5);\draw (11,1.3) ellipse (1 and 1.5);\end{scope}
\begin{scope}\clip (11,1.3) ellipse (1 and 1.5);\draw (11,-1.3) ellipse (1 and 1.5);\end{scope}

\begin{scope}\clip (12,-2) rectangle (14,2);\draw (13,0) ellipse (1.5 and 1);\end{scope}

\begin{scope}\clip (14,-2) rectangle ((17,2);\draw (15,0) ellipse (1.5 and 1);\end{scope}
\begin{scope}\clip (15,-.9) ellipse (1 and 1.5);\draw (15,1.3) ellipse (1 and 1.5);\end{scope}
\begin{scope}\clip (15,1.3) ellipse (1 and 1.5);\draw (15,-1.3) ellipse (1 and 1.5);\end{scope}\draw (13,0) node {$\cdots$};

\draw[-stealth] (8,0) -- (9,0);
\draw [decorate, decoration = {brace,mirror}] (1,-1.5) -- (7.5,-1.5);\draw (4,-1.5) node[below] {$$(g-1)/2$$};
\draw [decorate, decoration = {brace,mirror}] (9.5,-1.5) -- (16.5,-1.5);\draw (13,-1.5) node[below] {$$(g-1)/2$$};

\end{tikzpicture}
\\[.05 in]
\indent $e)$

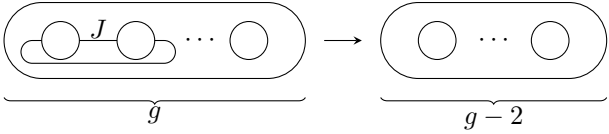
\begin{tikzpicture}[scale=1/2]
\draw (-1.5,1) -- (4.5,1) arc (90:-90:1) -- (-1.5,-1) arc (270:90:1);
\draw (-1,0) circle (.5); \draw (1,0) circle (.5); \draw (2.7,0) node {$\cdots$}; \draw (4,0) circle (.5);
\draw (-.5,0) -- (.5,0); \draw (1.5,0) -- (1.75,0) arc (90:-90:.3) -- (-1.75,-.6) arc (270:90:.3) -- (-1.5,0);
\node at (0,.3) {$J$};

\draw (8.5,1) -- (12.5,1) arc (90:-90:1) -- (8.5,-1) arc (270:90:1);
\draw (9,0) circle (.5); \draw (12,0) circle (.5); \draw (10.5,0) node {$\cdots$}; 

\draw[-stealth] (6,0) -- (7,0);
\draw [decorate, decoration = {brace,mirror}] (-2.5,-1.5) -- (5.5,-1.5);\draw (1.5,-1.5) node[below] {$g$};
\draw [decorate, decoration = {brace,mirror}] (7.5,-1.5) -- (13.5,-1.5);\draw (10.5,-1.5) node[below] {$g-2$};
\end{tikzpicture}
\\[.05 in]
\indent $f)$

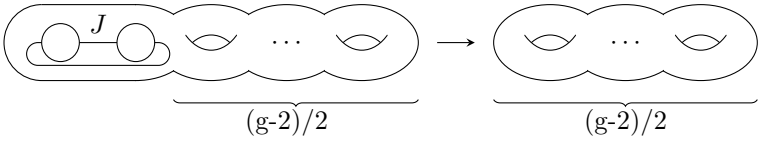
\begin{tikzpicture}[scale=1/2]
\begin{scope}\clip (0,-2) rectangle (1,2);\draw (0,0) ellipse (1.5 and 1);\end{scope}

\draw (-2.5,1) -- (0,1);\draw (0,-1) --(-2.5,-1) arc(270:90:1); 
\draw (-2,0) circle (.5); \draw (0,0) circle (.5);
\draw (-1.5,0) -- (-.5,0); \draw (.5,0) -- (.6,0) arc (90:-90:.3) -- (-2.6,-.6) arc (270:90:.3) -- (-2.5,0);
\node at (-1,.5) {$J$};

\begin{scope}\clip (1,-2) rectangle (3,2);\draw (2,0) ellipse (1.5 and 1);\end{scope}
\begin{scope}\clip (2,-.9) ellipse (1 and 1.5);\draw (2,1.3) ellipse (1 and 1.5);\end{scope}
\begin{scope}\clip (2,1.3) ellipse (1 and 1.5);\draw (2,-1.3) ellipse (1 and 1.5);\end{scope}\draw (4,0) node {$\cdots$};

\begin{scope}\clip (3,-2) rectangle (5,2);\draw (4,0) ellipse (1.5 and 1);\end{scope}

\begin{scope}\clip (5,-2) rectangle ((8,2);\draw (6,0) ellipse (1.5 and 1);\end{scope}
\begin{scope}\clip (6,-.9) ellipse (1 and 1.5);\draw (6,1.3) ellipse (1 and 1.5);\end{scope}
\begin{scope}\clip (6,1.3) ellipse (1 and 1.5);\draw (6,-1.3) ellipse (1 and 1.5);\end{scope}

\begin{scope}\clip (9.5,-2) rectangle (12,2);\draw (11,0) ellipse (1.5 and 1);\end{scope}
\begin{scope}\clip (11,-.9) ellipse (1 and 1.5);\draw (11,1.3) ellipse (1 and 1.5);\end{scope}
\begin{scope}\clip (11,1.3) ellipse (1 and 1.5);\draw (11,-1.3) ellipse (1 and 1.5);\end{scope}

\begin{scope}\clip (12,-2) rectangle (14,2);\draw (13,0) ellipse (1.5 and 1);\end{scope}

\begin{scope}\clip (14,-2) rectangle ((17,2);\draw (15,0) ellipse (1.5 and 1);\end{scope}
\begin{scope}\clip (15,-.9) ellipse (1 and 1.5);\draw (15,1.3) ellipse (1 and 1.5);\end{scope}
\begin{scope}\clip (15,1.3) ellipse (1 and 1.5);\draw (15,-1.3) ellipse (1 and 1.5);\end{scope}\draw (13,0) node {$\cdots$};

\draw[-stealth] (8,0) -- (9,0);
\draw [decorate, decoration = {brace,mirror}] (1,-1.5) -- (7.5,-1.5);\draw (4,-1.5) node[below] {$$(g-2)/2$$};
\draw [decorate, decoration = {brace,mirror}] (9.5,-1.5) -- (16.5,-1.5);\draw (13,-1.5) node[below] {$$(g-2)/2$$};

\end{tikzpicture}
\\[.05 in]
\indent $g)$

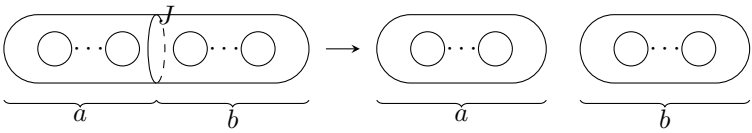
\begin{tikzpicture}[scale=.45]

\draw (-1.5,1) -- (5.5,1) arc (90:-90:1) -- (-1.5,-1) arc (270:90:1);
\draw (-1,0) circle (.5); \draw (1,0) circle (.5); \draw (0,0) node {$\cdots$}; \draw (3,0) circle (.5);\draw (4,0) node {$\cdots$}; \draw (5,0) circle (.5);
\draw [dashed] (2,1) arc (90:-90:0.25 and 1);
\draw (2,1) arc (90:270: 0.25 and 1); \node at (2.3,1) {$J$};

\draw (9.5,1) -- (12.5,1) arc (90:-90:1) -- (9.5,-1) arc (270:90:1);
\draw (10,0) circle (.5); \draw (12,0) circle (.5); \draw (11,0) node {$\cdots$}; 
\draw (15.5,1) -- (18.5,1) arc (90:-90:1) -- (15.5,-1) arc (270:90:1);
\draw (16,0) circle (.5); \draw (18,0) circle (.5); \draw (17,0) node {$\cdots$}; 

\draw[-stealth] (7,0) -- (8,0);
\draw [decorate, decoration = {brace,mirror}] (-2.5,-1.5) -- (2,-1.5);\draw (-.25,-1.5) node[below] {$a$};
\draw [decorate, decoration = {brace,mirror}] (2,-1.5) -- (6.5,-1.5);\draw (4.25,-1.5) node[below] {$b$};
\draw [decorate, decoration = {brace,mirror}] (8.5,-1.5) -- (13.5,-1.5);\draw (11,-1.5) node[below] {$a$};
\draw [decorate, decoration = {brace,mirror}] (14.5,-1.5) -- (19.5,-1.5);\draw (17,-1.5) node[below] {$b$};
\end{tikzpicture}
\\[.05 in]
\indent $h)$

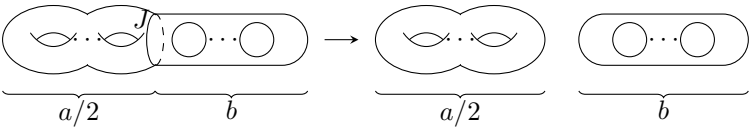
\begin{tikzpicture}[scale=.45]
\begin{scope}\clip (-2,-2) rectangle (1,2);\draw (0,0) ellipse (1.5 and 1);\end{scope}
\begin{scope}\clip (0,-.9) ellipse (1 and 1.5);\draw (0,1.3) ellipse (1 and 1.5);\end{scope}
\begin{scope}\clip (0,1.3) ellipse (1 and 1.5);\draw (0,-1.3) ellipse (1 and 1.5);\end{scope}
 
\begin{scope}\clip (1,-2) rectangle (3,2);\draw (2,0) ellipse (1.5 and 1);\end{scope}
\begin{scope}\clip (2,-.9) ellipse (1 and 1.5);\draw (2,1.3) ellipse (1 and 1.5);\end{scope}
\begin{scope}\clip (2,1.3) ellipse (1 and 1.5);\draw (2,-1.3) ellipse (1 and 1.5);\end{scope}\draw (1,0) node {$\cdots$};
\draw [dashed] (3,.75) arc (90:-90:0.25 and .75);
\draw (3,.75) arc (90:270: 0.25 and .75); \node at (2.6,.6) {$J$};

\draw (3,.75) -- (6.75,.75) arc (90:-90:.75) -- (3,-.75);
\draw (4,0) circle (.5); \draw (6,0) circle (.5); \draw (5,0) node {$\cdots$}; 

\begin{scope}\clip (9,-2) rectangle (12,2);\draw (11,0) ellipse (1.5 and 1);\end{scope}
\begin{scope}\clip (11,-.9) ellipse (1 and 1.5);\draw (11,1.3) ellipse (1 and 1.5);\end{scope}
\begin{scope}\clip (11,1.3) ellipse (1 and 1.5);\draw (11,-1.3) ellipse (1 and 1.5);\end{scope}
 
\begin{scope}\clip (12,-2) rectangle (15,2);\draw (13,0) ellipse (1.5 and 1);\end{scope}
\begin{scope}\clip (13,-.9) ellipse (1 and 1.5);\draw (13,1.3) ellipse (1 and 1.5);\end{scope}
\begin{scope}\clip (13,1.3) ellipse (1 and 1.5);\draw (13,-1.3) ellipse (1 and 1.5);\end{scope}\draw (12,0) node {$\cdots$};

\draw (16.25,.75) -- (19.75,.75) arc (90:-90:.75) -- (16.25,-.75)arc(270:90:.75);
\draw (17,0) circle (.5); \draw (19,0) circle (.5); \draw (18,0) node {$\cdots$}; 

\draw[-stealth] (8,0) -- (9,0);
\draw [decorate, decoration = {brace,mirror}] (-1.5,-1.5) -- (3,-1.5);\draw (.75,-1.5) node[below] {$a/2$};
\draw [decorate, decoration = {brace,mirror}] (3,-1.5) -- (7.5,-1.5);\draw (5.25,-1.5) node[below] {$b$};
\draw [decorate, decoration = {brace,mirror}] (9.5,-1.5) -- (14.5,-1.5);\draw (12,-1.5) node[below] {$a/2$};
\draw [decorate, decoration = {brace,mirror}] (15.5,-1.5) -- (20.5,-1.5);\draw (18,-1.5) node[below] {$b$};

\end{tikzpicture}
\\[.05 in]

\qed
\end{Proof}

\begin{Corollary}
Starting from $og$ or $ng$, the manifold decomposition game ends after at most $2g$ moves. 
The smallest number of moves is $k$ for $n(2k)$ and $k+1$ for $n(2k+1)$. 
\end{Corollary}
\begin{Proof}
All moves decrease the genus sum except cases $b$ and $g$. In these two cases, the genus sum stays the same and the number of 
components that are not spheres increases by $1$. Hence, the games with the most moves would use these two cases repeatedly for $g$ moves 
until there are $g$ components that are not spheres left, each of genus $1$. After $g$ more moves using either case $a$ or $c$ there will be $g$ 
surfaces left of genus $0$, i.e. spheres. There are no moves from a sphere, so the game ends. 

Starting from $n(2k)$, use case $e$ repeatedly to get a sphere after $k$ moves. Starting from $n(2k+1)$, use case $d$ once to get 
$ok$ and then use case $a$ repeatedly to get a sphere after $k$ more moves. 
\qed
\end{Proof}

We will next show that
starting from an orientable surface, the $G$-series is $01\dot{2}\dot{0}$ according to increasing genus, and
starting from a nonorientable surface, the $G$-series is $012\dot{4}60\dot{3}$ according to increasing genus. 
Nim addition determines the $G$-values of the remaining compact $2$-manifolds. The $G$-series is composed of Grundy numbers
which we call $G$-values or values for short. The dots over digits indicate a repeating block. 

\begin{Theorem}
The $G$-series for $og$ is $01\dot{2}\dot{0}$.
\end{Theorem}

\begin{Proof}
From Theorem 2.1 there are only two possible moves from $og$ for $g>0$: $og \rightarrow o(g-1)$, and $og \rightarrow (oa,ob)$ for $a,b>0, g=a+b.$
This game is equivalent to the octal game 4.3 in Guy and Smith~\cite{GS}. Their game starts with a heap of $g$ counters and in each move a player may 
split a heap into two heaps or remove one counter from a heap. For completeness, we provide a proof of the claim. 

The proof is by induction on $g$. Initial $G$-values are computed below. \\
$G(o0)=G(\emptyset) =0$, $G(o1)=mex\{G(o0)\}=mex\{0\}=1$,\\
 $G(o2)=mex\{G(o1), G(o1,o1)\} =mex\{1,G(o1)\oplus G(o1)\} =mex\{1,1\oplus1\}=mex\{1,0\}=2$.

Assume $g>2$ and $G$-values for smaller genus have been established. \\[.05 in]
Case 1: $g$ is even. Then $a,b$ are both even, or both odd. If both are even, then by hypothesis $G(a)$ and $G(b)$ are both 2, and their nim sum is 0. 
If both are odd then either one of $a$ or $b$ is 1, or both are larger than 1. Then the value is either the nim sum of 1 and 0 or the nim sum of 0 and 0.
Thus $G(g)=mex\{Go(g-1),0,1,0\}=mex\{0,0,1,0\}=2$.\\[.05 in]
Case 2: $g$ is odd. Then $a$ and $b$ have opposite parity with $G$-values 0 or 1, and 2, and nim sum 2 or 3. Thus $G(g)=mex\{Go(g-1),2,3\}=mex\{2,2,3\}=0.$
\qed
\end{Proof}

\begin{Theorem}
The $G$-series for $ng$ is $012\dot{4}60\dot{3}$.
\end{Theorem}

\begin{Proof}
Table 1 shows the computation of $G$-values of $ng$ for $g=0,\ldots,14$. For each $g$, each move from $ng$ is given along with the moves $G$-value. The $G$-value of 
$ng$, $G(ng)$, is the minimal excluded  non-negative integer, or mex, of these move values. We will show by induction, that for $g$ sufficiently large, 
the set of $G$-values for moves from $ng$ is the same as the set of $G$-values for moves from $n(g+4)$. It will follow that $G(ng)=G(n(g+4))$. 

The possible moves from $ng$ and $n(g+4)$ correspond to cases $c,\ldots,h$ in Theorem 2.1. We consider each case in turn.\\[.05 in]
Case c. $ng \rightarrow n(g-1)$ and $n(g+4) \rightarrow n(g+3)$. If the $G$-values of the right sides are equal, then the $G$-values of the left sides will be equal. The proof will proceed by induction. \\[.05 in]
Case e. $ng \rightarrow n(g-2)$ and $n(g+4) \rightarrow n(g+2)$. Similar to case c.\\[.05 in]
Case d. $ng \rightarrow o((g-1)/2)$  if $g-1 \ge 0$ and even, and $n(g+4) \rightarrow o((g+3)/2)$. Similar to above but using the previous theorem for the $G$-series for $o(g)$.\\[.05 in]
Case f. $ng \rightarrow o((g-2)/2)$  if $g-2 \ge 0$ and even, and $n(g+4) \rightarrow o((g+2)/2)$.  Similar to case d. \\[.05 in]
Case g. We list the moves in two tables, for $g=2k$ and for $g=2k+1.$ 

\begin{tabular}{lll}
$n(2k) \rightarrow $     & $(n1,n(2k-1)),$            & $\dots ,A=(n(k-2),n(k+2)),$ \\
                                     & $B=(n(k-1),n(k+1)),$    & $(nk,nk).$                            \\
$n(2k+4) \rightarrow $ & $(n1,n(2k+3)),$           & $\dots ,(n(k-2),n(k+6)), $       \\
                                     & $(n(k-1),n(k+5)),$        & $(nk,n(k+4)),$                      \\
                                     & $B^{\prime}=(n(k+1),n(k+3)),$   & $A^{\prime}=(n(k+2),n(k+2)).$         \\
\end{tabular}

The moves from $n(2k)$ correspond inductively to moves from $n(2k+4)$ in the order listed. The two additional moves from $n(2k+4)$ labeled $A^{\prime}$ and $B^{\prime}$ correspond to moves from $n(2k)$ labeled $A$ and $B$. 

\begin{tabular}{lll}
$n(2k+1) \rightarrow $     & $(n1,n(2k)),$            & $\dots ,(n(k-2),n(k+3)),$ \\
                                     & $C=(n(k-1),n(k+2)),$    & $D=(nk,n(k+1).$             \\
$n(2k+5) \rightarrow $ & $(n1,n(2k+4)),$           & $\dots ,(n(k-2),n(k+7)), $    \\
                                     & $(n(k-1),n(k+6)),$        & $(nk,n(k+5)),$                      \\
                                     & $D^{\prime}=(n(k+1),n(k+4)),$   & $C^{\prime}=(n(k+2),n(k+3)).$         \\
\end{tabular}

Same as above except $C$ corresponds to $C^{\prime}$ and $D$ to $D^{\prime}$. \\

\noindent Case h. We list the moves in two tables, for $g=2k$ and for $g=2k+1.$ 

\begin{tabular}{lll}
$n(2k) \rightarrow $     & $(o1,n(2k-2)),$        & $\dots ,A=(o(k-3),n6)),$   \\
                                     & $B=(o(k-2),n4),$    & $C=(o(k-1),n2).$                \\
$n(2k+4) \rightarrow $ & $(o1,n(2k+2)),$      & $\dots ,(o(k-3),n10), $       \\
                                     & $(o(k-2),n8),$         & $A^{\prime}=(o(k-1),n6),$                \\
                                     & $B^{\prime}=(ok,n4),$         & $C^{\prime}=(o(k+1),n2).$             \\
\end{tabular}

The moves from $n(2k)$ correspond inductively to moves from $n(2k+4)$ in the order listed, with the exception of the rule labeled $C$, as $G(n2) \neq G(6)$. 
The two additional moves from $n(2k+4)$ labeled $B^{\prime}$ and $C^{\prime}$ correspond to moves from $n(2k)$ labeled $B$ and $C$. Note that for
$A$ and $A^{\prime}$ to correspond we must have $k>4$ to use the periodicity of the $G$-series of $o(g)$. So $2k+4>12$.

\begin{tabular}{lll}
$n(2k+1) \rightarrow $     & $(o1,n(2k-1)),$     &  $\dots ,D=(o(k-2),n5)),$ \\
                                     &  $E=(o(k-1),n3)),$    &   $F=(ok,n1),$                   \\
$n(2k+5) \rightarrow $ & $(o1,n(2k+3)),$       & $\dots ,(o(k-1),n7), $       \\
                                     & $D^{\prime}=(ok,n5),$           & $E^{\prime}=(o(k+1),n3),$            \\
                                     & $F^{\prime}=(o(k+2),n1).$     &                                        \\
\end{tabular}

Same as above except $D$ corresponds to $D^{\prime}$, $E$ to $E^{\prime}$, and $F$ to $F^{\prime}$.

Table 1 shows the computation of $G$-values of $ng$ for $g=0,\ldots,14$. We see that $G(ng)=G(n(g+4))$ for $3 \leq g \leq 10$. 
Assume $g \ge 9$ and that for all $g^{\prime}$ with $12 \le g^{\prime}\le g+3$, the $G$-value of $ng^{\prime}$ is given.
Then for each case above, the set of $G$-values for moves from $ng$ is the same as the set of $G$-values for moves from $n(g+4)$. 
Thus, the set of $G$-values for all moves from $ng$ is the same as the set of $G$-values for all moves from $n(g+4)$, and so $G(ng)=G(n(g+4))$. 
\qed

\end{Proof}

\newpage

{\setlength\tabcolsep{1 pt} 
\begin{center}
\begin{longtable}{c|ccc|c}
\caption{\textbf{Calculating initial G-values of ng}}\label{tab:long} \\

\hline \multicolumn{1}{c|}{\textbf{g}} & \multicolumn{1}{c}{\textbf{$\cdots$}} & \multicolumn{1}{c}{\textbf{Moves from ng and G-values}} &  \multicolumn{1}{c|}{\textbf{$\cdots$}} & \multicolumn{1}{c}{\textbf{G(ng)}} \\ \hline 
\endfirsthead

\multicolumn{5}{c}%
{{\bfseries \tablename\ \thetable{} -- continued from previous page}} \\
\hline \multicolumn{1}{c|}{\textbf{g}} & \multicolumn{1}{c}{\textbf{$\cdots$}} & \multicolumn{1}{c}{\textbf{Moves from ng and G-values}} &  \multicolumn{1}{c|}{\textbf{$\cdots$}} & \multicolumn{1}{c}{\textbf{G(ng)}} \\ \hline 
\endhead

\hline \multicolumn{5}{|r|}{{Continued on next page}} \\ \hline
\endfoot

\hline \hline
\endlastfoot

0   &                     & $\emptyset$    &               &  0      \\
1   & n0, 0            &                    &                    &  1      \\
2   & n1,1             & n0, 0           & (n1, n1), 0   &  2      \\  \hline
3   & n2, 2            & n1, 1           & o1, 1          & 4       \\
          &  (n1, n2), 3   & (o1, n1), 0   &                  &               \\
4   & n3, 4            & n2, 2            & o1, 1           & 6       \\
          & (n1, n3), 5    & (n2, n2), 0    & (o1, n2), 3  &               \\
5   & n4, 6            & n3, 4            & o2, 2           & 0       \\
          & (n1, n4), 7   & (n2, n3), 6    & (o1, n3), 5   &               \\
          & (o2, n1), 3   &                     &                    &               \\
6   & n5, 0            & n4, 6            & o2, 2           &  3      \\
          & (n1, n5), 1     & (n2, n4), 4     & (n3, n3), 0   &                \\
          & (o1, n4), 7     &  (o2, n2), 0    &                   &                \\  \hline
7   & n6, 3            & n5, 0            & o3, 0           & 4      \\
          & (n1, n6), 2    & ( n2, n5), 2    & (n3, n4), 2  &                \\
          & (o1, n5), 1    & (o2, n3), 6    & ( o3, n1), 1  &                \\
8   & n7, 4            & n6, 3            & o3, 0           & 6        \\
         & (n1, n7), 5     &  (n2, n6), 1   & (n3, n5), 4   &                \\
         & (n4, n4), 0     &  (o1, n6), 2   & (o2, n4), 4   &                \\
         & (o3, n2), 2     &                   &                      &                \\
9   & n8, 6            & n7, 4            & o4, 2           &  0       \\
         & (n1, n8), 7     & (n2, n7), 6    & (n3, n6), 7   &                \\
         & (n4, n5), 6     & (o1, n7), 5    & (o2, n5), 2   &                \\
         & (o3, n3), 4     & (o4, n1), 3    &                   &             \\
10 & n9, 0            & n8, 6            & o4, 2          & 3       \\
          & (n1, n9), 1    & (n2, n8), 4    & (n3, n7), 0   &              \\
          & (n4, n6), 5   & (n5, n5), 0    & (o1, n8), 7   &              \\
          &(o2, n6), 1   & (o3, n4), 6    & (o4, n2), 0   &              \\ \hline
11  & n10, 3        & n9, 0            & o5, 0           & 4      \\
         & (n1, n10), 2 & (n2, n9), 2    & (n3, n8), 2   &            \\
           & (n4, n7), 2  & (n5, n6), 3   & (o1, n9), 1    &            \\
          & (o2, n7), 6 & (o3, n5), 0    & (o4, n3), 6   &              \\
          & (o5, n1), 1  &                   &                    &               \\
12 & n11, 4          & n10, 3          & o5, 0          &   6      \\
         &(n1, n11), 5   & (n2, n10), 1   &(n3, n9), 4    &               \\
         &(n4, n8), 0      &(n5, n7), 4    & (n6, n6), 0    &               \\
          &  (o1, n10), 2    &  (o2, n8), 4    & (o3, n6), 3   &           \\
          & (o4, n4), 4      &  (o5, n2), 2  &                       &          \\      
13  & n12, 6        & n11, 4            & o6, 2           & 0      \\
         & (n1, n12), 7 & (n2, n11), 6    & (n3, n10), 7   &            \\
           & (n4, n9), 6  & (n5, n8), 6   &  (n6 n7), 7   &            \\
          & (o1, n11), 5 & (o2, n9), 2    & (o3, n7), 4  &              \\
          &(o4, n5), 2   &  (o5, n3), 4  &   (o6, n1), 3   &               \\
14 & n13, 0          & n12, 6          & o6, 2          &   3      \\
         &(n1, n13), 1   & (n2, n12), 4   &(n3, n11), 0    &               \\
         &(n4, n10), 5      &(n5, n9), 0   & (n6, n8), 5    &               \\
         &(n7, n7), 0       &(o1, n12), 7  &(o2, n10), 1      &           \\
          &  (o3, n8), 6   &  (o4, n6), 1     & (o5, n4), 6   &           \\
          & (o6, n2), 0      &                     &                       &          \\      \hline

\end{longtable}
\end{center}
}

\end{document}